\documentclass[12pt]{amsart}
\usepackage{amsmath,amsfonts,amsthm,amscd,amssymb,euscript,
graphics,psfrag}

\newcommand{\beq}{\begin{equation}}

\newcommand{\eeq}{\end{equation}}

\def\mmtrx#1#2#3#4{\begin{bmatrix} #1 & #2  \\ #3 & #4\end{bmatrix}}
\def\mtrx#1#2#3#4#5#6#7#8#9{\begin{bmatrix} #1 & #2 & #3\\ #4 & #5 & #6 \\ #7 & #8 & #9\end{bmatrix}}

 \theoremstyle{definition}

\newtheorem{definition}{Definition}

 \newtheorem{conjecture} {Conjecture}
\newtheorem*{utheorem}{Theorem}

\theoremstyle{plain}      

 \newtheorem{proposition} {Proposition}
 \newtheorem{theorem} {Theorem}

\newcommand{\natls}{{\mathbb N}}
\newcommand{\col}{{\hbox{col}}}
\newcommand{\row}{{\hbox{row}}}
\newcommand{\vol}{{\hbox{vol}}}

\newcommand{\Tr}{{\hbox{\rm Tr}}}

\newcommand{\Sc}{{\hbox{\rm Sc}}}

\newcommand{\reals}{{\mathbb R}}
\newcommand{\zed}{{\mathbb Z}}
\newcommand{\E}{{\mathbb E}}

\newcommand{\Eu}{\E_{U_N}}

\newcommand{\card}{{\rm card}}

\newcommand{\Sp}{{\rm Sp}}

\title[Pseudomoments of the Riemann  zeta-function]%
{Pseudomoments of the Riemann  zeta-function and pseudomagic squares}

\author{Brian Conrey}
\address{American Institute of Mathematics
360 Portage Avenue, Palo Alto, CA 94306} \email{conrey@aimath.org}

\author{Alex Gamburd}
\address {Department of
Mathematics, University of California, Santa Cruz and  Department
of Mathematics, Stanford University, Stanford, CA 94305 }
\email{agamburd@math.stanford.edu}

\thanks{The first author was supported in part by an NSF grant.
The second author was supported in part by the NSF postdoctoral
fellowship. }

\begin{document}

\begin{abstract}
We compute integral moments of partial sums of the Riemann zeta
function on the critical line and obtain an expression for the
leading coefficient as a product of the standard arithmetic factor
and a geometric factor. The geometric factor is equal to the
volume of the convex polytope of substochastic matrices  and is
equal  to the leading coefficient in the expression for moments of
truncated characteristic polynomial of a random unitary matrix.

\end{abstract}
\maketitle

\section{Introduction}
\label{sec:intro}
\subsection{Moments of the Riemann zeta function}

The Riemann zeta-function is defined  for $Re(s)
>1$ by
\beq
\zeta(s)=\sum_{n=1}^{\infty}\frac{1}{n^s}=\prod_{p}\left(1-\frac{1}{p^s}\right)^{-1}.
\eeq As is well-known \cite{Ti},  $\zeta$ has meromorphic
continuation to the whole complex plane with a single simple pole
at $s=1$ with residue $1$.  Further, it satisfies a functional
equation, relating the value of $\zeta(s)$ and the value of
$\zeta(1-s)$, \beq \zeta(s)=\chi(s) \zeta(1-s), \eeq where \beq
\chi(s)=2^{s} \pi^{s-1} \sin{\frac{\pi s}{2}} \Gamma(1-s). \eeq
Following the standard notation we write $s=\sigma + i t.$

The problem of computing the moments of $\zeta$  on the critical
line $\sigma=\frac{1}{2}$  is fundamental, difficult and
longstanding.

The second moment was obtained by Hardy and Littlewood \cite{HL}
in 1918: \beq \frac{1}{T}\int_{0}^{T}\left|\zeta\left(\frac{1}{2}+
it \right)\right|^2 \, dt \sim \log {T} , \eeq the fourth moment
was obtained by Ingham \cite{I26} in 1926: \beq
\frac{1}{T}\int_{0}^{T}\left|\zeta\left(\frac{1}{2}+ it
\right)\right|^4 \, dt \sim \frac{1}{2 \pi^{2}}\log^{4}{T} .\eeq
The asymptotics of higher moments is not known.  It has long been
conjectured that \beq
\frac{1}{T}\int_{0}^{T}\left|\zeta\left(\frac{1}{2}+ it\right)
\right|^{2k} \, dt \sim c_k \log^{k^2}{T}. \eeq In 1984 Conrey and
Ghosh \cite{CG84} gave the moment conjecture a more precise form;
namely, they conjectured  that there should be a factorization
\beq \label{ggg} c_k = \frac{g_k a_k}{\Gamma(1+k^2)}, \eeq where
$a_k$ is an arithmetic factor given by \beq  \label
{e:af}a_k=\prod_{p}\left(1-\frac{1}{p}\right)^{k^2}
\sum_{j=0}^{\infty}\frac{d_k(p^j)^2}{p^j}, \eeq and $g_k$ is a
geometric factor, which should be an integer. Using Dirichlet
polynomial techniques Conrey and Ghosh \cite{CG98} conjectured
that $g_3=42$ and Conrey and Gonek \cite{CoGo01} conjectured that
$g_4=24024.$

\subsection{The Riemann zeta-function and characteristic polynomials
of random matrices}

In the past few years, following the work of Keating and Snaith
\cite{KS00a}, Conrey and Farmer \cite{CF00}, Hughes, Keating and
O'Connell \cite{HKO1, HKO2}, and Conrey, Farmer, Keating,
Rubinstein, and Snaith \cite{CFKRS1} it has become clear that the
leading order asymptotic of the moments of the Riemann zeta
function can be conjecturally understood in terms of corresponding
quantities of the characteristic polynomial of the random unitary
matrices.
 Let $M$ be a matrix in $U(N)$ chosen uniformly with respect to
 Haar measure, denote by $e^{i \theta_1}, \dots, e^{i \theta_N}$
 its eigenvalues, and consider the characteristic polynomial of
 $M$:
 \beq \label{cpu}
 P_M(z)=\det(M-zI)=\prod_{j=1}^{N}(e^{i \theta_j}-z).
\eeq

Keating and Snaith (see also Baker and Forrester \cite{BF97})
computed the moments of $P_M$ with respect to Haar measure on
$U(N)$ and found that \beq
 M_N(s)=E_{U(N)}|P_M(z)|^{2s}=
\prod_{j=1}^{N}\frac{\Gamma(j)\Gamma(j+2s)}{\Gamma(j+s)^2}. \eeq
They also showed that \beq \lim_{N \to
\infty}\frac{M_N(s)}{N^{s^2}} = \frac{G(1+s)^2}{G(1+2s)}, \eeq
where $G(s)$ is Barnes double Gamma function satisfying $G(1)=1$
and $G(z+1)=\Gamma(z) G(z).$ For $s=k$ an integer \beq
\frac{G(1+k)^2}{G(1+2k)}=\prod_{j=0}^{k-1}\frac{j!}{(j+n)!} .\eeq
For $k=1, 2, 3$ the quantity above is in agreement with the value
of $g_k$ in the  theorems of Hardy and Littlewood, and Ingham and
the conjecture of Conrey and Ghosh.

The conjecture of Keating and Snaith \cite{KS00a} (considerably
refined and extended in \cite{CFKRS1}) is as follows:
\begin{conjecture} [Keating and Snaith \cite{KS00a}]  \label{con:ks}
 \beq \label{e:ks}
\frac{1}{T}\int_{0}^{T}\left|\zeta\left(\frac{1}{2}+ it\right)
\right|^{2k} \, dt \sim a_k g_k\log^{k^2}{T}, \eeq where $a_k$ is
an arithmetic factor given by \eqref{e:af} and $g_k$ is a
``geometric'' factor (here the notation is different from
\eqref{ggg}) given by \beq g_k  \label{e:gk} =\lim_{N\to\infty}
\frac{E_{U(N)}|P_M(z)|^{2k}}{N^{k^2}}=\prod_{j=0}^{k-1}\frac{j!}{(j+n)!}.\eeq

\end{conjecture}

 \subsection{Characteristic polynomials of unitary matrices and
 magic squares}
 The moments of the secular coefficients of the random unitary
 matrices have also been recently investigated.  If $M$ is
 a random unitary matrix, following the notation preceding
 equation \eqref{cpu} we write:
 \beq \label{cpu2}
 P_M(z)=\det(M-zI)=\prod_{j=1}^{N}(e^{i \theta_j}-z)=
 (-1)^N \sum_{j=0}^{N}\Sc_j(M)z^{N-j}(-1)^j,
 \eeq
where $\Sc_j(M)$ is the $j$-th \emph{secular coefficient} of the
characteristic polynomial. Note that \beq \Sc_1(M)=\Tr(M),\eeq and
\beq \Sc_N(M)=\det(M).\eeq

Moments of the higher secular coefficients were studied by Haake
and collaborators \cite{Ha1, Ha2}  who obtained: \beq
\label{e:haake} \E_{U(N)} \Sc_j(M) =0,  \quad \E_{U(N)}
|\Sc_j(M)|^2 =1; \eeq and posed the question of computing the
higher moments. The answer is given by Theorem \ref{thm:agpd},
which we state below after pausing  to give the following
definition.
\begin{definition} If $A$ is an $m$ by $n$ matrix with nonnegative
integer entries and with row and column sums
$$r_i=\sum_{j=1}^{n}a_{ij},$$
$$c_j=\sum_{i=1}^{m}a_{ij};$$
then the the row-sum vector $\row(A)$ and column-sum vector
$\col(A)$ are defined by
$$\row(A)=(r_1, \dots, r_m),$$
$$\col(A)=(c_1, \dots, c_n).$$
\end{definition}
Given two partitions $\mu=(\mu_1, \dots, \mu_m)$ and
$\tilde{\mu}=(\tilde{\mu}_1, \dots, \tilde{\mu}_n)$ (see
\cite{Mac} for the partition notation)  we denote by $N_{\mu
\tilde{\mu}}$ the number of nonnegative integer matrices $A$ with
$\row(A) = \mu$ and $\col(A) = \tilde{\mu}$.

For example, for $\mu=(2, 1, 1)$ and $\tilde{\mu}=(3, 1)$  we have
$N_{\mu \tilde{\mu}}=3;$ and the matrices in question are
$$\begin{bmatrix}
2&0\\
1&0\\
0&1
\end{bmatrix},
\begin{bmatrix}
2&0\\
0&1\\
1&0
\end{bmatrix},
\begin{bmatrix}
1&1\\
1&0\\
1&0
\end{bmatrix}.
$$

For $\mu=(2, 2, 1)$ and $\tilde{\mu}=(3, 1,1)$  we have $N_{\mu
\tilde{\mu}}=8;$ and the matrices in question are
$$\begin{bmatrix}
0&1&1\\
2&0&0\\
1&0&0
\end{bmatrix},
\begin{bmatrix}
1&1&0\\
1&0&1\\
1&0&0
\end{bmatrix},
\begin{bmatrix}
1&0&1\\
1&1&0\\
1&0&0
\end{bmatrix},
\begin{bmatrix}
2&0&0\\
0&1&1\\
1&0&0
\end{bmatrix},$$
$$\begin{bmatrix}
2&0&0\\
1&1&0\\
0&0&1
\end{bmatrix},
\begin{bmatrix}
2&0&0\\
1&0&1\\
0&1&0
\end{bmatrix},
\begin{bmatrix}
1&1&0\\
2&0&0\\
0&0&1
\end{bmatrix},
\begin{bmatrix}
1&0&1\\
2&0&0\\
0&1&0
\end{bmatrix}.$$

We are ready to state the following Theorem, proved in \cite{DiGa}

\begin{theorem}[\cite{DiGa}] \label{thm:agpd}
(a)  Consider $\mathbf{a}=(a_1, \dots, a_l)$ and $\mathbf{b}=(b_1,
\dots, b_l)$ with $a_j$, $b_j$ nonnegative natural numbers. Then
for $N \geq \max\left(\sum_1^l ja_j, \sum_1^l j b_j \right)$ we
have \beq \label{e:mixedmom} \Eu \prod\limits^l_{j=1}
(\Sc_j(M))^{a_j} \overline{(\Sc_j(M))} ^{b_j} = N_{\mu
\tilde{\mu}}. \eeq Here  $\mu$ and $\tilde{\mu}$ are partitions
$\mu=\langle 1^{a_1}\dots l^{a_l}\rangle $, $\tilde{\mu}=\langle
1^{b_1}\dots l^{b_l}\rangle $  and $N_{\mu \tilde{\mu}}$ is the
number of nonnegative integer matrices $A$ with $\row(A) = \mu$
and $\col(A) = \tilde{\mu}$.

 (b)
In particular, for $N \geq j k$ we have \footnote{We remark that
in \cite{fg} the answer is also obtained in the case $N < j k$: it
is related to enumeration of magic squares with certain additional
constraints.} \beq \label{cuemom} E_{U(N)}|\Sc_j(M)|^{2k}= H_k(j),
\eeq where $H_k(j)$ is the number of $k\times k$ nonnegative
integer matrices with each row and column summing up to $j$ --
``magic squares''.
\end{theorem}

\subsection{Magic Squares}
The reader is likely to have encountered objects, which following
Ehrhart~\cite{EE77} are referred to as ``historical  magic
squares''. These are square matrices of order $k$,  whose entries
are nonnegative integers ($1, \dots, k^2$) and whose rows and
columns sum up to the same number.  The oldest such object,
\begin{equation}  \label{hm}
\mtrx{4}{9}{2}{3}{5}{7}{8}{1}{6}
\end{equation}
first appeared in ancient Chinese literature under the name
\emph{Lo Shu} in the third millennium BC and repeatedly reappeared
in the cabbalistic and occult literature in the middle ages. Not
knowing ancient Chinese, Latin, or Hebrew it is difficult to
understand what is ``magic'' about \emph{Lo Shu}; it is quite easy
to understand however why it keeps reappearing: there is (modulo
reflections) only one historic magic square of order $3$.

Following MacMahon \cite{PM15} and Stanley \cite{RS73}, what is
referred to as  magic squares in modern combinatorics are square
matrices of order $k$,  whose entries are nonnegative integers and
whose rows and columns sum up to the same number $j$. The number
of magic squares of order $k$ with row and column sum $j$, denoted
by $H_k(j)$, is of great interest; see \cite{DG} and references
therein. The first few values are easily obtained:
\begin{equation} \label{hk1} H_k(1)=k!, \eeq corresponding to all
 $k$ by $k$ permutation matrices (this is the $k$-th moment of the
 traces of powers leading in the work of Diaconis and Shahshahani
 \cite{DS94}
 to the result on the asymptotic normality);
\beq \label{e:h1}
 H_1(j)=1,
\eeq corresponding to $1 \times 1$ matrix $[j]$ (this is the
result of Haake and collaborators given in equation
\eqref{e:haake}).  We also easily obtain  $ H_2(j)= j+1, $
corresponding to $ \mmtrx{i}{j-i}{j-i}{i} $, but  the value of
$H_3(j)$ is considerably more involved:
\begin{equation}\label{h3j}
H_3(j)= \binom{j+2}{4} +\binom{j+3}{4} +\binom{j+4}{4}. \eeq This
expression was first obtained by Mac Mahon in 1915 and a simple
proof was found only a few years ago by M. Bona. The main results
on $H_k(j)$ are given by the following theorems, proved by Stanley
and Ehrhart (see \cite{EE73, EE77, RS73, RS74, stanley86}):
\begin{utheorem}[Stanley] \label{stan}
 $H_k(j)$ is a polynomial in $j$ of degree $(k-1)^2,$  having
 ``trivial zeroes'' at the negative integers,
\beq \label{e:rec1} H_k(-1)=H_k(-2)=\dots=H_k(-k+1)=0,  \eeq
 and satisfying the following ``functional equation'':

  \beq \label{e:rec2}
H_k(-k-j)=(-1)^{k-1}H_k(j).\eeq

\end{utheorem}

It can be shown that the statements above are equivalent to \beq
\sum_{j\geq 0}H_k(j) x^{j}=\frac{h_0+h_1 x +\dots +h_d
x^{d}}{(1-x)^{(k-1)^2 +1}}, \quad d=k^2-3k+2, \eeq with
$h_0+h_1+\dots h_d \ne 0$ and $h_i=h_{d-i}$.

For example,
$$
H_3(j)=\frac{1}{8}j^4+\frac{3}{4}j^3+\frac{15}{8}j^2+\frac{9}{4}j+1.
$$
and
$$
\sum_{j\geq 0}H_3(j) x^{j}=\frac{1+x+x^2}{(1-x)^{5}}.
$$

$$
\sum_{j\geq 0}H_4(j) x^{j}=\frac{1+14 x +87x^2 +148x^3+87 x^4+14
x^5 +x^6}{(1-x)^{10}}.
$$

\begin{utheorem}[Ehrhart]
The leading coefficient of $H_k(j)$ is the relative volume of
$\mathcal{B}_{k}$ - the $k$-th Birkhoff polytope, i.e. leading
coefficient is equal to $\frac{\vol(\mathcal{B}_{k})}{k^{k-1}}$.
\end{utheorem}
By definition, the $k$-th Birkhoff polytope is the convex hull of
permutation matrices: \beq \mathcal{B}_k= \left\{ (x_{ij}) \in
\reals^{k^2} \biggm | x_{ij}\geq 0; \quad \sum_{i=1}^{k}x_{ij} =
1; \quad \sum_{j=1}^{k}x_{ij} = 1 \right\} .\eeq
 In the example above,
$\vol(\mathcal{B}_{3})=\frac{1}{8}\times 9$.

\subsection{Pseudomoments of the Riemann zeta-function and
pseudomagic squares}

 The purpose of this paper is to prove the
following result:
\begin{theorem}  \label{thm:main} Let $a_k$ be  the arithmetic factor given by equation
\eqref{e:af}.
  Then \beq \label{e:main} \lim_{T \to \infty}
\frac{1}{T}\int_{0}^{T} \left \vert \sum_{n=1}^{X}
\frac{1}{n^{\frac{1}{2}+i t}} \right \vert^{2k} d t= a_k \gamma_k
(\log X)^{k^2} + O\left( (\log X)^{k^2-1} \right). \eeq

Here
$\gamma_k$ is the geometric factor,  $\gamma_k
=\vol(\mathcal{P}_k)$, where $\mathcal{P}_k$ is the convex
polytope in $\reals^{k^2}$ defined by the following inequalities:

\beq   \label{e:pkt} \mathcal{P}_k= \left\{ (x_{ij}) \in
\reals^{k^2} \biggm | x_{ij}\geq 0; \quad \sum_{i=1}^{k}x_{ij}
\leq 1; \quad \sum_{j=1}^{k}x_{ij} \leq 1 \right\} .\eeq
\end{theorem}

The connection with the characteristic polynomials of unitary
matrices is as follows.  From Theorem \ref{thm:agpd} it follows
that if we consider \emph{truncated} characteristic polynomial

\beq \label{tpu}
 P_{M,l}(z)=\sum_{j=0}^{l}\Sc_j(M)z^{N-j}(-1)^j,
 \eeq
we have for $N \ge lk$ \beq \label{e:tpu1} \E_{U(N)}
|P_{M,l}(z)|^{2k} =G_k(l), \eeq where $G_k(l)$ denotes the number
of $k \times k$ nonnegative integer matrices with row and column
sums less than or equal to $l$ (referred to as  ``pseudomagic
squares'' by Ehrhart \cite{EE77}) : \beq \label{e:gkt} G_k(l)
=\card \left\{ (x_{ij}) \in \zed^{k^2} \biggm | x_{ij}\geq 0;
\quad \sum_{i=1}^{k}x_{ij} \leq l; \quad \sum_{j=1}^{k}x_{ij} \leq
l \right\} .\eeq

Ehrhart \cite{EE77}  proved that $G_k(l)$ is a
 polynomial in $l$ of degree $k^2$ with leading coefficient
 given by $\gamma_k=\vol(\mathcal{P}_k)$; in fact $G_k(l) =\card \left(l \mathcal{P}_k
\cap \zed^{k^2}\right).$   For example,
$$
G_2(l)=\frac{1}{6}(l+1)(l+2)(l^2+3l+3),
$$
and
$$\vol(\mathcal{P}_2)=\frac{1}{6}.
$$

Hence we can rewrite the geometric factor $\gamma_k$ in a manner
similar to the expression for $g_k$ in \eqref{e:gk} as follows:
\beq \label{e:gammak} \gamma_k= \lim_{l\to\infty}
\frac{E_{U(lk)}|P_{M,l}(z)|^{2k}}{l^{k^2}}.\eeq

The proof proceeds as follows.  In section \ref{sec:ps} we obtain
an expression for $\gamma_k$ in terms of a multiple complex
integral. In section \ref{sec:proof} we express the left-hand side
of \eqref{e:main} as a multiple complex integral and then show
that the leading terms in the two resulting expressions are equal.

\noindent \textbf{Acknowledgements.}  It is a pleasure to thank
Dan Bump for pointing out that \eqref{e:tpu1} is a consequence of
Theorem \ref{thm:agpd}.

\section{Pseudomagic squares} \label{sec:ps}
Let $G_k(l)$ denote the number of $k \times k$ nonnegative integer
matrices with row and column sums less than or equal to $l$ given
by \eqref{e:gkt}
 (we remark that $H_{k+1}(l)$, the number of magic squares, is
 obtained by imposing an additional diophantine inequality
 $\sum_{i, j} x_{ij} \geq (k-1)l.$)

We have the following expression for $G_k(l)$ as a multiple
complex integral:

\begin{proposition}\label{p:op}
  Notation being as above we have
\beq \label{e:gklexpr1} G_k(l)=\frac{1}{(2 \pi i)^{2k}}
\idotsint_{\substack{|w_i|=\epsilon_i\\|z_j|=\epsilon_j}}
\frac{(w_1 \dots w_k z_1 \dots z_k)^{-l-1} \prod_{i=1}^{k} dw_{i}
\prod_{j=1}^{k} dz_{j}}{\prod_{i, j}(1-w_i z_j) \prod_{i=1}^{k}
(1-w_{i})\prod_{j=1}^{k}(1-z_{j})}. \eeq
\end{proposition}

The proof follows the approach in \cite{MB00}, which we now
review.

Let $\zed^n$ denote an $n$-dimensional integer lattice in
$\reals^n$  and let $\mathcal{P}$ be a convex polytope in
$\reals^n$ whose vertices are on the lattice $\zed^n$
 ($\mathcal{P}_k$ is a convex lattice polytope in
$\zed^{k^2}$).

 Any convex lattice polytope situated in the nonnegative orthant can be described as an
intersection of finitely many half-spaces: \beq \label{e:polp}
\mathcal{P} = \left\{\mathbf{x} \in \reals_{\geq 0}^n \bigm |
\mathbf{A} \mathbf{x} \leq \mathbf{b} \right \}, \eeq where
$\mathbf{A}$ is an $m \times n$ integer matrix and $\mathbf{b} \in
\zed^{m}$. Consider now the function of an integer-valued variable
$l$ describing  the number of lattice points that lie inside the
dilated polytope $l \mathcal{P}$: \beq L(\mathcal{P}, l) = \card
\{l \mathcal{P} \cap \zed^n \};\eeq with this notation $G_k(l) =
L(\mathcal{P}_k, l).$ Denote the columns of $\mathbf{A}$ by
$\mathbf{c}_1, \dots, \mathbf{c}_n$.   Using multivariate
generating functions it is proved in \cite{MB00} that for the
lattice polytope $\mathcal{P}$ given by \eqref{e:polp} we have the
following expression for $L(\mathcal{P}, l)$: \beq
\label{e:mulint1} L(\mathcal{P}, l)=\frac{1}{(2 \pi i)^{m}}
\idotsint_{|z_j|=\epsilon_j} \frac{\prod_{j=1}^{m}z_j^{-l
b_j-1}}{\prod_{l=1}^{n}(1-\mathbf{z}^{\mathbf{c}_{l}})
\prod_{j=1}^{m}(1-z_j)}  \, d \mathbf{z}. \eeq In the expression
above we use the standard multivariate notation
$\mathbf{x}^{\mathbf{y}}=x_1^{y_1}\dots x_n^{y_n}$.

Now for $\mathcal{P}_k$  the defining system of diophantine
inequalities is given in \eqref{e:pkt}; the corresponding
$\mathbf{A}$ is a $(2k \times k^{2})$ matrix given by \beq
\mathbf{A}= \left(
\begin{matrix}
1 & \dots & 1 & & & & & & &\\
& &  &1 &\dots  &1  & & & &\\
&  &  & & & &\ddots & & &\\
&  &  & & & & &1 &\dots &1\\
1&  &  &1 & & & &1 & &\\
&\ddots  &  & &\ddots & &\dots & &\ddots &\\
&  &1 & & &1 & & & &1
\end{matrix} \right)
, \eeq and $\mathbf{b} =(1, \dots, 1) \in \zed^{2k}.$ Proposition
\ref{p:op} now follows from \eqref{e:mulint1}; for notational
convenience we have split the variables into two groups $w_1,
\dots, w_k$ and $z_1, \dots, z_k$.


\section{Proof of the Theorem} \label{sec:proof}
By the mean-value theorem for Dirichlet polynomials due to
Montgomery and Vaughan \cite{MV} we have \beq \label{grs} \lim_{T
\to \infty} \frac{1}{T}\int_{0}^{T} \left \vert \sum_{n=1}^{X}
\frac{1}{n^{\frac{1}{2}+i t}} \right \vert^{2k}\, dt
=\sum_{n=1}^{X}\frac{d_{k, X}^{2}(n)}{n}, \eeq where $d_{k, X}(n)$
is defined by:
$$d_{k, X}(n)=\sum_{\substack{l_1 \dots l_k = n\\
l_1, \dots, l_k \leq X}}1 .$$

Consequently we have:

\beq \label{e:mulsum} \lim_{T \to \infty} \frac{1}{T}\int_{0}^{T}
\left \vert \sum_{n=1}^{X} \frac{1}{n^{\frac{1}{2}+i t}} \right
\vert^{2k}\, dt =
\sum_{\substack{1\leq l_i \leq X\\
1 \leq m_j \leq X \\l_1 \dots l_k =m_1 \dots
m_k}}\frac{1}{\sqrt{l_1 \dots l_k m_1 \dots m_k}}. \eeq

Now we use the discontinuous integral \beq \label{e:int}\frac{1}{2
\pi i}\int_{c-i\infty}^{c+i \infty}\frac{X^s}{s} \, d s =
\begin{cases} 0 , &\text{if $0 < X <1$;}\\
1, &\text{if $X > 1$}
\end{cases}
\eeq where $c>0$ to pick the terms of the Dirichlet series which
are less than $X$.  Denoting the integral in equation
\eqref{e:int} by $\int\limits_{(c)}$ we can now express the
right-hand side of \eqref{e:mulsum} as follows: \beq
\label{e:mulint}\frac{1}{(2 \pi i)^{2k}}\int\limits_{(2)} \dots
\int\limits_{(2)} \prod_{i=1}^{k}\frac{X^{u_i}}{u_i}
\prod_{j=1}^{k}\frac{X^{v_j}}{v_j} F(u_1, \dots, u_k, v_1, \dots,
v_k) \, d u_1\dots d u_{k} d v_1\dots d v_{k}, \eeq where \beq
F(u_1, \dots, u_k, v_1, \dots, v_k)=
\sum_{\substack{l_i \geq 1\\
m_j \geq  1 \\
l_1 \dots l_k =m_1 \dots m_k}}\frac{1}{l_1^{\frac{1}{2}+u_1} \dots
l_k^{\frac{1}{2}+u_k} m_1^{\frac{1}{2}+v_1} \dots
m_k^{\frac{1}{2}+v_k}}. \eeq To simplify notation let
$\mathbf{u}=(u_1, \dots, u_k)$, $\mathbf{v}=(v_1, \dots, v_k)$,
$d\mathbf{u}=d u_1\dots d u_{k}$ and $d\mathbf{v}=d v_1\dots d
v_{k}.$

Now since for a multiplicative function $g(n)$ we have the Euler
product identity:\beq
\sum_{n=1}^{\infty}g(n)=\prod_{p}(1+g(p)+g(p^2)+g(p^3)+\dots),
\eeq it follows that
\begin{multline} \label{e:tozeta}
F(u_1, \dots, u_k, v_1, \dots, v_k)=
\prod_{p}\left(\sum_{n=1}^{\infty}\sum_{\substack{\alpha_1+\dots+\alpha_k=n\\
\beta_1+\dots+\beta_k=n\\\alpha_i\geq 0, \, \beta_j\geq
0}}\frac{1}{p^{\alpha_1(\frac{1}{2}+u_1)+\dots+\beta_k(\frac{1}{2}+v_k)}}\right)\\
=\prod_{p}\left(1+\sum_{i,
j}\frac{1}{p^{1+u_i+v_j}}+\dots\right)=G(\mathbf{u},
\mathbf{v})\prod_{i, j}\zeta(1+u_i+v_j),
\end{multline}

where $G$ is an Euler product which is absolutely convergent for
$|u_i|< \frac{1}{4}$, $|v_j|< \frac{1}{4}$.

Since
$$\sum_{\substack{\alpha_1+\dots+\alpha_k=n\\
\beta_1+\dots+\beta_k=n\\\alpha_i\geq 0, \, \beta_j\geq
0}}1=d_k^2(p^n),$$ if we let all $u_i$ and $v_j$ be equal to
$\delta$ we obtain \beq G(\delta, \dots
\delta)=\prod_{p}(1-\frac{1}{p^{2\delta+1}})^{k^2}
\sum_{n=0}^{\infty}d_k^2(p^n)p^{-2 n\delta - n}, \eeq and
consequently \beq \label{e:guva} \lim_{\mathbf{u}, \mathbf{v} \to
0}G(\mathbf{u}, \mathbf{v})=
\prod_{p}\left(1-\frac{1}{p}\right)^{k^2}
\sum_{n=0}^{\infty}\frac{d_k(p^n)^2}{p^n}=a_k.\eeq

Summarizing, we have obtained the following expression for the
left-hand side of \eqref{e:main}:
\begin{multline} \label{e:summ} \lim_{T
\to \infty} \frac{1}{T}\int_{0}^{T} \left \vert \sum_{n=1}^{X}
\frac{1}{n^{\frac{1}{2}+i t}} \right \vert^{2k}\, dt=
\\
\frac{1}{(2 \pi i)^{2k}}\int \limits_{c} \dots \int \limits_{c}
G(\mathbf{u}, \mathbf{v})\prod_{i, j}\zeta(1+u_i+v_j)
\frac{X^{\sum(u_i +v_j)}}{ \prod_{i, j}u_i v_j} \, d \mathbf{u} \,
d \mathbf{v}.
\end{multline}
Now using the fact that $(s-1) \zeta(s)$ is analytic in the entire
complex plane together with the standard techniques and bounds
pertaining to $\zeta$, we obtain that the leading term in
\eqref{e:summ} is given by \beq \label{e:mta} \frac{a_k}{(2 \pi
i)^{2k}}\int \limits_{c} \dots \int\limits_{c} \frac{X^{\sum(u_i
+v_j)}}{ \prod_{i, j}(1-e^{-u_i-v_j}) \prod_{i, j}u_i v_j} \, d
\mathbf{u} \, d \mathbf{v}, \eeq where we have used
\eqref{e:guva}.

Write \beq \label{e:expan} \frac{1}{\prod_{i,
j}(1-e^{-u_i-v_j})}=\prod_{i, j}\left[\sum_{a_{ij} \geq 0}
(e^{-u_i-v_j})^{a_{ij}}\right]. \eeq A term $e^{-\mathbf{u}
\mathbf{\alpha}} e^{-\mathbf{v} \mathbf{\beta}}$ in this expansion
is obtained by choosing an $\natls$-matrix $A^{t} =(a_{ij})^{t}$
of finite support with $\row(A) = \mathbf{\alpha}$ and
$\col(A)=\mathbf{\beta}$. Hence the coefficient of $e^{-\mathbf{u}
\mathbf{\alpha}} e^{-\mathbf{v} \mathbf{\beta}}$ in
\eqref{e:expan} is the number $N_{\mathbf{\alpha}\mathbf{\beta}}$
of $\natls$-matrices $A$ with $\row(A) = \mathbf{\alpha}$ and
$\col(A)=\mathbf{\beta}$: \beq \label{e:expan2} \frac{1}{\prod_{i,
j}(1-e^{-u_i-v_j})}
=\sum_{\mathbf{\alpha}\mathbf{\beta}}N_{\mathbf{\alpha}\mathbf{\beta}}
e^{-\mathbf{u} \mathbf{\alpha}} e^{-\mathbf{v} \mathbf{\beta}}.
\eeq

Further,  let $l = \log X$  and rewrite the integral \eqref{e:int}
as follows: \beq \label{e:int2}\frac{1}{2 \pi
i}\int_{c-i\infty}^{c+i \infty}\frac{e^{ls}}{s} \, d s =
\begin{cases} 0 , &\text{if $l < 0$;}\\
1, &\text{if $l > 0$}.
\end{cases}
\eeq

We now express the integral appearing in \eqref{e:mta} using
\eqref{e:expan2} and apply \eqref{e:int2} to obtain:
\begin{multline} \label{e:mta2} \frac{1}{(2 \pi i)^{2k}}\int\limits_{c} \dots
\int\limits_{c}\ \frac{X^{\sum(u_i +v_j)}}{ \prod_{i,
j}(1-e^{-u_i-v_j}) \prod_{i, j}u_i v_j} \, d \mathbf{u} \, d
\mathbf{v} =
\\\frac{1}{(2 \pi i)^{2k}}\int\limits_{c} \dots
\int\limits_{c} \prod_{i}\frac{e^{l u_i}\, d
u_i}{u_i}\prod_{j}\frac{e^{l v_j}\, d
v_j}{v_j}\sum_{\mathbf{\alpha}\mathbf{\beta}}N_{\mathbf{\alpha}\mathbf{\beta}}
e^{-\mathbf{u} \mathbf{\alpha}} e^{-\mathbf{v} \mathbf{\beta}}=\\
 \sum_{\substack{\mathbf{\alpha} \leq \mathbf{l}\\
\mathbf{\beta}\leq \mathbf{l}}} N_{\mathbf{\alpha}\mathbf{\beta}}=
\card \left\{ (x_{ij}) \in \natls^{k^2} \biggm | \quad
\sum_{i=1}^{k}x_{ij} \leq l; \quad \sum_{j=1}^{k}x_{ij} \leq l
\right\} = G_k(l).
\end{multline}

We remark that this proves that the integrals given by
\eqref{e:mta2} and \eqref{e:gklexpr1} are equal; a direct proof
using, for example, a change of variables has thus far eluded us.
We also remark that the integral expression for $G_k(l)$ given by
\eqref{e:gklexpr1}  has served only as a motivation for the proof
presented above.  We also note that sums related to the expression
given by the right-hand side of \eqref{grs} were considered in
\cite{GS00}.

\section{Generalizations}  Note that in fact we have proved
\beq \label{e:main1} \lim_{T \to \infty} \frac{1}{T}\int_{0}^{T}
\left \vert \sum_{n=1}^{X} \frac{1}{n^{\frac{1}{2}+i t}} \right
\vert^{2k} d t= a_k G_k(\log X) + O\left( (\log X)^{k^2-1}
\right). \eeq The proof given in the previous section easily
generalizes to yield the following result:

\begin{theorem}  Let $a_k$ be  the arithmetic factor given by equation
\eqref{e:af}.
  Then up to the lower order terms we have: \beq \label{e:main2} \lim_{T \to \infty}
\frac{1}{T}\int_{0}^{T} \left \vert \sum_{n=1}^{X_1}
\frac{1}{n^{\frac{1}{2}+i t}} \right \vert^{2} \dots \left \vert
\sum_{n=1}^{X_k} \frac{1}{n^{\frac{1}{2}+i t}} \right \vert^{2} d
t \sim  a_k G_k(\log{X_1}, \dots, \log{X_k}). \eeq

Here we assume that $X_i=Y^{m_i}$ with $m_i = O(1)$ and $Y \to
\infty$;   \beq G_k(l_1, \dots, l_k) = \card \left\{ (x_{ij}) \in
\natls^{k^2} \biggm | \quad \sum_{i=1}^{k}x_{ij} \leq l_j; \quad
\sum_{j=1}^{k}x_{ij} \leq l_i \right\} .\eeq
\end{theorem}

Finally, we note that in \cite{DiGa} results analogous to Theorem
\ref{thm:agpd} are proved for orthogonal and symplectic group; for
example the result for symplectic group is as follows:

\begin{theorem} \label{thm:spms}
(a)  Consider $\mathbf{a}=(a_1, \dots, a_l)$ with $a_j$
nonnegative natural numbers. Let  $\mu$  be  a partition
$\mu=\langle 1^{a_1}\dots l^{a_l}\rangle.$  Then for $N \geq
\sum_1^l ja_j$ and $|\mu|$ \emph{even} we have \beq
\label{e:mixedmomsp} \E_{\Sp(2N)}\prod\limits^l_{j=1}
(\Sc_j(M))^{a_j} = NSP_{\mu}. \eeq
  Here $NSP_{\mu}$ is the number of nonnegative
\emph{symmetric} integer matrices $A$ with $\row(A) = \col(A) =
\mu$ and  with all diagonal entries of $A$ even.

 (b)
In particular, for $N \geq j k$  and $jk$ even we have \beq
\label{csemom} E_{\Sp(2N)}\Sc_j(M)^{k}= S^{sp}_k(j), \eeq where
$S^{sp}_k(j)$ is the number of $k\times k$ symmetric nonnegative
integer matrices with each row and column summing up to $j$  and
all diagonal entries even (equivalently, the number of  j-regular
 graphs on k vertices with loops and multiple edges).
\end{theorem}

 We will
present  analogues of Theorem \ref{thm:main} for L-functions with
orthogonal and symplectic symmetries in a forthcoming paper.  Here
we state a representative result for $L(s,\chi_d)$ with
$\chi_d(n)=(\tfrac{d}{n})$ where $d$ is  a fundamental
discriminant, which has symplectic
 symmetry.

 \begin{theorem} \label{thm:lsymp}
Let $b_k$ be  the arithmetic factor given by
$$
b_k=\prod_p \frac{(1-\tfrac 1p)^{\frac{k(k+1)}{2}}} {1+\tfrac 1p}
\left(\frac{(1-\tfrac{1}{\sqrt{p}})^{-k}+
(1+\tfrac{1}{\sqrt{p}})^{-k}}{2}+\frac{1}{p}\right)
$$

  Then \beq \label{e:msymp} \lim_{T \to \infty}
\frac{1}{T^{*}} \sum_{d<T} \left (\sum_{n<X}
\frac{\chi_{d}(n)}{\sqrt{n}} \right )^{k} = \frac{6}{\pi^2}b_k F_k
(\log X) +O(\log X^{k^2+k-2/2}) .\eeq  Here
$F_k(l)$ is the polynomial in $l$ of degree $k(k+1)/2$ equal to
the  the number of $k \times k$ symmetric nonnegative integer
matrices with row and column sums less than or equal to $l$ and
all diagonal entries even.
\end{theorem}

The connection with the characteristic polynomials of symplectic
matrices is as follows.  From Theorem \ref{thm:spms} it follows
that if we consider \emph{truncated} characteristic polynomial

$$
 P_{M,l}(z)=\sum_{j=0}^{l}\Sc_j(M)z^{N-j}(-1)^j,
$$
we have for $N \ge lk$  $$\E_{\Sp(2N)} P_{M,l}(z)^{k} =F_k(l);
$$ from results of
Ehrhart \cite{EE77} it follows that $F_k(l)$ is a
 polynomial in $l$ of degree $k(k+1)/2.$

\end{document}